\documentclass[12pt]{article}

\usepackage[cp866]{inputenc}
\usepackage{citehack}
\usepackage{amsmath}
\usepackage{amsfonts}
\usepackage{amssymb}
\newtheorem{theorem}{Theorem}

\newtheorem{ex}{Example}

\def\i{{\hat{i}}}
\def\j{{\hat{j}}}
\def\k{{\hat{k}}}

\def\ii{{\hat{\hat{i}}}}
\def\jj{{\hat{\hat{j}}}}
\def\kk{{\hat{\hat{k}}}}

\def\zr{\ltimes}

\def\Real{\mathbb{R}}
\def\Co{\mathbb{C}}

\def\g{\mathfrak{g}}
\def\h{\mathfrak{h}}

\def\so{\mathfrak{so}}
\def\spin{\mathfrak{spin}}
\def\simil{\mathfrak{sim}}

\def\gl{\mathfrak{gl}}
\def\su{\mathfrak{su}}
\def\u{\mathfrak{u}}
\def\sp{\mathfrak{sp}}
\def\f{\mathfrak{f}}

\def\z{\mathfrak{z}}
\def\R{\mathcal{R}}
\def\P{\mathcal{P}}

\def\id{\mathop\text{\rm id}\nolimits}
\def\spa{\mathop\text{{\rm span}}\nolimits}

\def\tr{\mathop\text{\rm tr}\nolimits}

\def\pr{\mathop\text{\rm pr}\nolimits}

\def\Ric{\mathop{{\rm Ric}}\nolimits}

\def\tRic{\mathop{\widetilde{\rm Ric}}\nolimits}

\def\p{\partial}

\begin{document}

\title{Holonomy of Einstein Lorentzian manifolds}
\author{Anton S. Galaev}


\maketitle

\centerline{Dedicated to Mark Volfovich Losik on his 75th birthday}

\begin{abstract}
  The classification of all possible holonomy algebras of Einstein and vacuum Einstein Lorentzian manifolds is obtained. 
It is shown that each such algebra
appears as the holonomy algebra of an Einstein (resp., vacuum Einstein) Lorentzian manifold, the direct constructions are
given. Also  the holonomy algebras of totally Ricci-isotropic Lorentzian manifolds are classified.
The classification of the  holonomy algebras of Lorentzian manifolds is reviewed and  a complete description of the spaces of
curvature tensors for these holonomies is given.
\end{abstract}


\section{Introduction} In contrast to the case of Riemannian manifolds, where the classification  of holonomy algebras 
 is a classical result, which has a lot of consequences and  applications both in geometry and physics (e.g. Riemannian manifolds with most of holonomy algebras are automatically Einstein or vacuum Einstein), see
\cite{Al,Ber,Besse,Jo00,Jo07} and the references therein, the classification of the holonomy algebras of Lorentzian
manifolds is achieved only recently \cite{BB-I,Leistner,Gal5,ESI} (we recall it in Section \ref{holclas}).  The most interesting case is
when a Lorentzian manifold $(M,g)$ admits a parallel distribution of isotropic lines and the manifold is locally indecomposable, i.e.
locally it is not a product of a Lorentzian and a Riemannian manifold. In this case the holonomy algebra is contained in
the maximal subalgebra $\simil(n)=(\Real\oplus\so(n))\zr\Real^n$ of the Lorentzian algebra $\so(1,n+1)$ preserving an
isotropic line (the dimension of $M$ is $n+2$). There is a number of recent physics literature dealing with these manifolds, see e.g. \cite{BCH,BCH1,CGHP,CFH,CHPP,FF00,G-P,Gibbons09,LS5}. In particular, in \cite{BCH,BCH1,Gibbons09} expressed the hope that Lorentzian manifolds with the holonomy algebras contained in $\simil(n)$ will find many applications in physics, e.g. they are of interest in M-theory and string theory.

The fundamental equation of General Relativity is  the Einstein equation.
In the absence of matter it has the form 
 \begin{equation} \Ric=\Lambda g,\end{equation}
where $g$ is a Lorentzian metric on a manifold $M$, $\Ric$ is the Ricci tensor of the metric $g$ and $\Lambda \in \Real$ is the cosmological constant.
If the metric of  a Lorentzian  manifold $(M,g)$ satisfies this equation, then $(M,g)$ is called {\it an Einstein manifold}. If moreover $\Lambda=0$, then it is called {\it vacuum Einstein or Ricci-flat}.
In dimension 4 the solutions of this equation are obtained in \cite{GT,KG1,KG2}.
In \cite{Bryant,FF00} are found the holonomy algebras of vacuum Einstein Lorentzian spin manifolds up to dimension 11 admitting parallel spinors
and some methods of construction of such manifolds are introduced. In dimension 11 these manifolds are purely gravitational supersymmetric solutions of 11-dimensional supergravity. Other constructions are provided in \cite{LS5}.
Recently in \cite{G-P} the Einstein equation on Lorentzian manifolds with holonomy algebras contained in $\simil(n)$ is
studied and in some cases solutions are obtained. See \cite{G-P} for the discussion of the importance of this problem and the references there.

In the present paper we classify all possible holonomy algebras of Einstein and vacuum Einstein Lorentzian manifolds.
First, in \cite{ESI} it is proved that {\it if $(M,g)$ is Einstein, then its holonomy algebra coincides  either with $(\Real\oplus\h)\zr\Real^n\subset\simil(n)$, or with $\h\zr\Real^n\subset\simil(n)$.} Here $\h\subset\so(n)$ is the holonomy algebra of a Riemannian manifold. 
In general for such $\h\subset\so(n)$ there is an orthogonal decomposition
\begin{equation}\label{LM0A}\Real^{n}=\Real^{n_1}\oplus\cdots\oplus\Real^{n_s}\oplus\Real^{n_{s+1}}\end{equation} and the
corresponding decomposition into the direct sum of ideals
\begin{equation}\label{LM0B}\h=\h_1\oplus\cdots\oplus\h_s\oplus\{0\}\end{equation} such that $\h$ annihilates
$\Real^{n_{s+1}}$, $\h_i(\Real^{n_j})=0$ for $i\neq j$, and $\h_i\subset\so(n_i)$ is an irreducible subalgebra for $1\leq
i\leq s$. Moreover, the Lie subalgebras $\h_i\subset\so(n_i)$ are the holonomy algebras of Riemannian manifolds. The classification of the Riemannian holonomy algebras shows that each $\h_i\subset\so(n_i)$
is either one of the Lie algebras $\so(n_i)$, $\u(\frac{n_i}{2})$, $\su(\frac{n_i}{2})$, $\sp(\frac{n_i}{4})$,
$\sp(\frac{n_i}{4})\oplus\sp(1)$, $G_2\subset\so(7)$, $\spin(7)\subset\so(8)$, or it is a symmetric Berger algebra 
(the last are the Riemannian holonomy algebras such that any manifold with such holonomy is locally symmetric, these Lie algebras are listed e.g. in \cite{Besse}). In Section \ref{Einstein} we prove the following two theorems

{\bf Theorem} \ref{thRic0}. {\it If $(M,g)$  is vacuum Einstein, then one of the following holds:
\begin{description} \item[(1)] The holonomy algebra of $(M,g)$ coincides with $(\Real\oplus\h)\zr\Real^n$, and in the decomposition \eqref{LM0B}
of  $\h\subset\so(n)$ at least one subalgebra $\h_i\subset\so(n_i)$ coincides with one of the Lie
algebras $\so(n_i)$, $\u(\frac{n_i}{2})$, $\sp(\frac{n_i}{4})\oplus\sp(1)$ or with a symmetric Berger algebra.
\item[(2)] The holonomy algebra of $(M,g)$ coincides with $\h\zr\Real^n$, and in the decomposition \eqref{LM0B}
of $\h\subset\so(n)$ each subalgebra $\h_i\subset\so(n_i)$ coincides with one of the Lie algebras
$\so(n_i)$, $\su(\frac{n_i}{2})$, $\sp(\frac{n_i}{4})$, $G_2\subset\so(7)$, $\spin(7)\subset\so(8)$. \end{description}
}

{\bf Theorem} \ref{thEinstein}. {\it
 If $(M,g)$ is Einstein and not vacuum Einstein, then the holonomy algebra of $(M,g)$ 
coincides with $(\Real\oplus\h)\zr\Real^n$, and in 
the decomposition \eqref{LM0B} of  $\h\subset\so(n)$ each subalgebras $\h_i\subset\so(n_i)$ coincides
with one of the Lie algebras $\so(n_i)$, $\u(\frac{n_i}{2})$, $\sp(\frac{n_i}{4})\oplus\sp(1)$ or with a symmetric Berger algebra.
Moreover,  in the decomposition \eqref{LM0A} it holds $n_{s+1}=0$.
}

These theorems recover the previous results for manifolds of dimension 4 obtained in \cite{Hall-Lonie00,Schell}.

In Section \ref{secExamples} we construct local Einstein  metrics having holonomy algebras from Theorems \ref{thRic0} and \ref{thEinstein}.
We consider metrics written in Walker coordinates of the form
\begin{equation}\label{Walker1}
g=2dx^0dx^{n+1}+\sum^{n}_{i,j=1}h_{ij}(x^1,...,x^{n})dx^idx^j+f(x^0,...,x^{n+1})(dx^{n+1})^2,\end{equation} 
and
\begin{equation}\label{Walker2}
g=2dx^0dx^{n+1}+\sum^{n}_{i=1}(dx^i)^2+2\sum^{n}_{i=1}u^i(x^1,...,x^{n+1}) dx^i
dx^{n+1}+f(x^0,...,x^{n+1})(dx^{n+1})^2.\end{equation} Here $h=\sum^{n}_{i,j=1}h_{ij}(x^1,...,x^{n})dx^idx^j$
is a Riemannian metric. For the metric \eqref{Walker1} we choose $h$ to be an Einstein Riemannian metric with a cosmological constant $\Lambda\neq 0$
(resp., $\Lambda= 0$) and holonomy algebra $\h\subset\so(n)$ as in Theorem \ref{thEinstein} (resp., (2) of Theorem \ref{thRic0}). Choosing an appropriate function $f$, we get that the metric $g$ is Einstein with the cosmological constant $\Lambda$ and the holonomy algebra $(\Real\oplus\h)\zr\Real^n$ (resp., $\h\zr\Real^n$). Next, for a given holonomy algebra $\h\subset\so(n)$ as in  (1) (resp., (2)) of Theorem \ref{thRic0} we give an algorithm how to choose the functions $u^i$ and $f$ in \eqref{Walker2} to make the metric $g$ vacuum Einstein with the holonomy algebra $(\Real\oplus\h)\zr\Real^n$ (resp., $\h\zr\Real^n$). As examples, we provide vacuum Einstein Lorentzian metrics with the holonomy algebras $G_2\zr\Real^7\subset\so(1,8)$ and $\spin(7)\zr\Real^8\subset\so(1,9)$.

If $(M,g)$ is a spin Lorentzian manifold and it admits a parallel spinor, then it is totally Ricci-isotropic (but not necessary vacuum Einstein, unlike in the Riemannian case), i.e. the image of its Ricci operator is isotropic \cite{Bryant,FF00}. We prove the following two theorems 

{\bf Theorem} \ref{thRicisotr}. {\it
If $(M,g)$  is totally Ricci-isotropic, then its holonomy algebra is the same  as in Theorem \ref{thRic0}.}

{\bf Theorem} \ref{thRicisotr1}. {\it
 If the holonomy algebra of $(M,g)$ is $\h\zr\Real^n$ and in the decomposition \eqref{LM0B}
of $\h\subset\so(n)$ each subalgebra $\h_i\subset\so(n_i)$ coincides with one of the Lie algebras
 $\su(\frac{n_i}{2})$, $\sp(\frac{n_i}{4})$, $G_2\subset\so(7)$, $\spin(7)\subset\so(8)$, then $(M,g)$ is totally Ricci-isotropic.}

Recall that an indecomposable Riemannian manifold with the holonomy algebra different from $\so(n)$ and $\u(\frac{n}{2})$ is automatically Einstein or vacuum Einstein \cite{Besse}. This is not the case for Lorentzian manifolds. Indeed, using results of Section \ref{secExamples} it is easy to construct non-Einstein metrics having all holonomy algebras as in Theorems \ref{thRic0} and \ref{thEinstein}. On the other hand, Theorem \ref{thRicisotr1} shows that Lorentzian manifolds with some holonomy algebras are automatically totally Ricci-isotropic.

To prove the above theorems  we need the complete description of the space of curvature tensors for Lorentzian holonomy algebras,
i.e. the space of values of the curvature tensor of a Lorentzian manifold, we provide it in Section \ref{curv}. The study of these spaces is begun in \cite{Gal1} and finished recently in \cite{GalP}.

Necessary facts from the holonomy theory can be found e.g. in \cite{Besse,BCH1,ESI,Jo00,Jo07}. Remark that by a Riemannian (resp., Lorentzian) manifold
we understand a manifold with a field $g$ of positive definite (resp., of signature $(-,+,...,+)$) symmetric bilinear forms on the tangent spaces.

{\bf Acknowledgments.} I am thankful to D.\,V.\,Alekseevsky and Thomas Leistner for useful and stimulating discussions on
the topic of this paper. 

 I wish to express my gratitude to Mark Volfovich Losik, my first teacher of Differential Geometry, for his guidance, mentoring and careful attention to my study and work since 1998 when I was at the very first course at Saratov State University. 

 The work was supported by the grant  201/09/P039 of the
Grant Agency of Czech Republic   and by the grant MSM~0021622409 of the Czech Ministry of Education.

\section{Classification of the Lorentzian holonomy algebras}\label{holclas}

Let $(\Real^{1,n+1},\eta)$ be the Minkowski space of dimension $n+2$, where  $\eta$ is a metric on $\Real^{n+2}$ of
signature $(1,n+1)$. We consider $(\Real^{1,n+1},\eta)$ as the tangent space $(T_xM,g_x)$ to a Lorentzian manifold $(M,g)$
at a point $x$. We fix a basis $p,e_1,...,e_n,q$ of $\Real^{1,n+1}$ such that the only non-zero values of $\eta$ are
$\eta(p,q)=\eta(q,p)=1$ and $\eta(e_i,e_i)=1$.  We will denote by $\Real^n\subset\Real^{1,n+1}$ the Euclidean subspace
spanned by the vectors $e_1,...,e_n$.

Recall that a subalgebra $\g\subset \so(1,n+1)$ is called  {\it irreducible} if it does not preserve any proper subspace
of $\Real^{1,n+1}$; $\g$ is called  {\it weakly-irreducible} if it does not  preserve any non-degenerate proper subspace
of $\Real^{1,n+1}$. Obviously, if $\g\subset \so(1,n+1)$ is irreducible, then it is weakly-irreducible.
 From the Wu  Theorem \cite{Wu} it follows that any Lorentzian manifold $(M,g)$ is either locally a product of
the manifold $(\Real,-(dt)^2)$ and of a Riemannian manifold, or of a Lorentzian manifold with weakly-irreducible holonomy
algebra and of a Riemannian manifold. The Riemannian manifold can be further decomposed into the product of a flat
Riemannian manifold and of Riemannian manifolds with irreducible holonomy algebras. If the manifold $(M,g)$ is simply
connected and geodesically complete, then these decompositions are global. This allow us to consider locally
indecomposable Lorentzian manifolds, i.e. manifolds with weakly-irreducible holonomy algebras. For example, a locally
decomposable Lorentzian manifold  $(M,g)$ is  Einstein if and only if  locally it is a product of Einstein manifolds with
the same cosmological constants as $(M,g)$. The only irreducible Lorentzian holonomy algebra is the whole Lie algebra
$\so(1,n+1)$ (\cite{Ber}), so we consider weakly-irreducible not irreducible holonomy algebras.

Denote  by $\simil(n)$ the subalgebra of $\so(1,n+1)$ that preserves the isotropic line $\Real p$. The Lie algebra
$\simil(n)$ can be identified with the following matrix algebra: \begin{equation}\label{matsim}\simil(n)=\left\{\left. \left (\begin{array}{ccc} a
&X^t & 0\\ 0 & A &-X \\ 0 & 0 & -a \\
\end{array}\right)\right|\, a\in \Real,\, X\in \Real^n,\,A \in \so(n) \right\} .\end{equation}
The above matrix can be identified with the triple $(a,A,X)$.  We get the decomposition
$$\simil(n)=(\Real\oplus\so(n))\zr\Real^n,$$ which means that $\Real\oplus\so(n)\subset\simil(n)$ is a subalgebra and
$\Real^n\subset\simil(n)$ is an ideal, and the Lie brackets of $\Real\oplus\so(n)$ with $\Real^n$ are given by the
standard representation of $\Real\oplus\so(n)$ in $\Real^n$.
 We see that $\simil(n)$ is isomorphic to the Lie algebra of the Lie group of
similarity transformations of $\Real^n$. The explicit isomorphism on the group level is constructed in \cite{Gal2}.

 If a weakly-irreducible subalgebra $\g\subset \so(1,n+1)$ preserves a
degenerate proper subspace $U\subset \Real^{1,n+1}$, then it preserves the isotropic line $U\cap U^\bot$, and $\g$ is
conjugated to a weakly-irreducible subalgebra of $\simil(n)$. Let $\h\subset\so(n)$ be a subalgebra. Recall that $\h$ is a
compact Lie algebra and we have the decomposition $\h=\h'\oplus\z(\h)$, where $\h'$ is the commutant of $\h$ and $\z(\h)$
is the center of $\h$.

The next theorem gives the classification of weakly-irreducible not irreducible holonomy algebras of Lorentzian manifolds.

\begin{theorem}
A subalgebra $\g\subset \simil(n)$ is the weakly-irreducible holonomy algebra of a Lorentzian manifold if and only if it
is conjugated to one of the following subalgebras:
\begin{description}
\item[type 1.] $\g^{1,\h}=(\Real\oplus\h)\zr\Real^n$, where $\h\subset\so(n)$ is the holonomy algebra of a Riemannian manifold;

\item[type 2.] $\g^{2,\h}=\h\zr\Real^n,$ where $\h\subset\so(n)$ is the holonomy algebra of a Riemannian manifold;

\item[type 3.] $\g^{3,\h,\varphi}=\{(\varphi(A),A,0)|A\in\h\}\zr\Real^n,$ where $\h\subset\so(n)$ is
the holonomy algebra of a Riemannian manifold
 with $\z(\h)\neq\{0\}$, and  $\varphi :\h\to\Real$ is a non-zero linear map with
$\varphi|_{\h'}=0$;

\item[type 4.] $\g^{4,\h,m,\psi}=\{(0,A,X+\psi(A))|A\in\h,X\in\Real^m\},$ where  $0<m<n$ is an integer,
$\h\subset\so(m)$ is the holonomy algebra of a Riemannian manifold with $\dim\z(\h)\geq n-m$, a decomposition $\Real^n=\Real^m\oplus \Real^{n-m}$ is fixed, and $\psi:\h\to \Real^{n-m}$
is a surjective linear map with $\psi|_{\h'}=0$. \end{description}\end{theorem}

The subalgebra $\h\subset\so(n)$ associated to a weakly-irreducible Lorentzian holonomy algebra $\g\subset \simil(n)$  is
called {\it the orthogonal part} of $\g$.
Recall that for the subalgebra $\h\subset\so(n)$ there are the decompositions \eqref{LM0A} and \eqref{LM0B}. 
This theorem is the result of the papers \cite{BB-I,Leistner,Gal5}, see \cite{ESI} for the whole history.

\section{The spaces of curvature tensors}\label{curv}

Let $W$ be a vector space and $\f\subset \gl(W)$ a subalgebra. The vector  space $$\R(\f)=\{R\in \Lambda^2
W^*\otimes\f|R(u,v)w+R(v, w)u+R(w, u)v=0 \text{ for all } u,v,w\in W\}$$ is called {\it the space of curvature tensors of
type} $\f$.  A subalgebra $\f\subset\gl(W)$ is called {\it a Berger algebra} if $$\g=\spa\{R(u, v)|R\in\R(\f),\,u,v\in
W\},$$ i.e. $\g$ is spanned by the images of the elements $R\in\R(\f)$.

If there is a pseudo-Euclidean metric $\eta$ on $W$ such that $\f\subset\so(W)$, then any $R\in\R(\f)$ satisfies
\begin{equation}\label{symR} \eta(R(u,v)z,w)=\eta(R(z,w)u,v)\end{equation}
for all $u,v,z,w\in W$. We identify $\so(W)$ with $\Lambda^2 W$ such that it holds $(u\wedge v)(z)=\eta(u,z)v-\eta(v,z)u$
for all $u,v,z\in W$. E.g. under this identification the matrix from \eqref{matsim} corresponds to the bivector $-ap\wedge q+A-p\wedge X$, $A\in\so(n)\simeq\Lambda^2\Real^n$.  Equality \eqref{symR} shows that the map $R:\Lambda^2 W\to\f\subset\Lambda^2 W$ is symmetric with
respect to the metric on $\Lambda^2 W$. In particular, $R$ is zero on the orthogonal complement $\f^\perp\subset \Lambda^2
W$. Thus, $R\in\odot^2 \f$.

In \cite{Gal1} the following theorem is proved.

\begin{theorem}\label{thR}\cite{Gal1}
It holds: \begin{description}
\item[(1)]
each $R\in\R(\g^{1,\h})$ is uniquely given by $$\lambda\in\mathbb{R},\ v\in \mathbb{R}^n,\ P\in\P(\mathfrak{h}),\
R_0\in\mathcal{R}(\mathfrak{h}),\text { and }T\in End(\mathbb{R}^n)\text{ with }T^*=T$$ in the following way
\begin{align*}
R(p,q)=&(\lambda,0,v),\qquad R(x,y)=(0,R_0(x,y),P(y)x-P(x)y),\\ R(x,q)=&(\eta(v,x),P(x),T(x)),\qquad R(p,x)=0 \end{align*}
 for all $x,y\in\Real^n$;
\item[(2)]  $R\in\R(\g^{2,\h})$ if and only if $R\in\R(\g^{1,\h})$, $\lambda=0$ and $v=0$;
\item[(3)]  $R\in\R(\g^{3,\h,\varphi})$ if and only if
$R\in\R(\g^{1,\h})$, $\lambda=0$, $R_0\in\R(\ker\varphi)$, and  $\eta(x,v)=\varphi( P(x))$ for all $x\in\Real^n$;
\item[(4)]  $R\in\R(\g^{4,\h,m,\psi})$ if and only if $R\in\R(\g^{1,\h})$, $\lambda=0$,
 $v=0$, $R_0\in\R(\ker\psi)$, and $\pr_{\mathbb{R}^{n-m}}\circ T=\psi\circ P$.
\end{description}
\end{theorem}

Remark that the decomposition $\R(\g^{1,\h})=\Real\oplus\Real^n\oplus\odot^2\Real^n\oplus\R(\h)\oplus\P(\h)$ is
$\Real\oplus\h$-invariant, but not $\g^{1,\h}$-invariant. 
Recall that for the subalgebra $\h\subset\so(n)$ there are the decompositions \eqref{LM0A} and \eqref{LM0B}. 
In addition  we have the decompositions
$$\P(\h)=\P(\h_1)\oplus\cdots\oplus\P(\h_s)$$ and $$\R(\h)=\R(\h_1)\oplus\cdots\oplus\R(\h_s).$$ 

The spaces $\R(\h)$ for the holonomy algebras of Riemannian manifolds $\h\subset\so(n)$ are computed by D.~V.~Alekseevsky
in \cite{Al}. Let $\h\subset\so(n)$ be an irreducible subalgebra.
 The space  $\R(\h)$ admits the following decomposition into
$\h$-modules \begin{equation}\label{razlR}\R(\h)=\R_0(\h)\oplus\R_1(\h)\oplus\R'(\h),\end{equation} where $\R_0(\h)$
consists of the curvature tensors with zero Ricci tensor, $\R_1(\h)$ consists of tensors annihilated by $\h$ (this
space is zero or one-dimensional), $\R'(\h)$ is the complement to these two spaces. Any element of $\R'(\h)$ has zero scalar curvature and non-zero Ricci tensor. If $\R(\h)=\R_1(\h)$, then any
Riemannian manifold with the holonomy algebra $\h$ is locally symmetric (such $\h\subset\so(n)$ is called {\it a symmetric Berger algebra}). Note that any indecomposable locally symmetric Riemannian
manifold is Einstein and not vacuum Einstein. Note that $\R(\h)=\R_0(\h)$ if $\h$ is any of the algebras: $\su(\frac{n}{2})$,
$\sp(\frac{n}{4})$, $G_2\subset\so(7)$, $\spin(7)\subset\so(8)$. This implies that each Riemannian manifold with any of
these holonomy algebras is vacuum Einstein. Next,
$\R(\u(\frac{n}{2}))=\Real\oplus\R'(\u(\frac{n}{2}))\oplus\R(\su(\frac{n}{2}))$ and
$\R(\sp(\frac{n}{4})\oplus\sp(1))=\Real\oplus\R(\sp(\frac{n}{4}))$. Hence any  Riemannian manifold with the holonomy
algebra $\sp(\frac{n}{4})\oplus\sp(1)$ is Einstein and not vacuum Einstein, and a Riemannian manifold with the holonomy algebra
$\u(\frac{n}{2})$ can not be vacuum Einstein. Finally, if an indecomposable $n$-dimensional Riemannian manifold is vacuum Einstein,
then its holonomy algebra is one of $\so(n)$, $\su(\frac{n}{2})$, $\sp(\frac{n}{4})$,  $G_2\subset\so(7)$,
$\spin(7)\subset\so(8)$.

Now we turn to the space $\P(\h)$, where $\h\subset\so(n)$ is an irreducible subalgebra. Consider the $\h$-equivariant map
$$\tRic:\P(\h)\to\Real^n, \qquad \tRic(P)=\sum_{i=1}^nP(e_i)e_i.$$ This definition does not depend on the choice of the
orthogonal basis $e_1,...,e_n$ of $\Real^n$. Denote by $\P_0(\h)$ the kernel of $\tRic$ and let $\P_1(\h)$ be its
orthogonal complement in $\P(\h)$. Thus, $$\P(\h)=\P_0(\h)\oplus\P_1(\h).$$ Since $\h\subset\so(n)$ is irreducible and the
map $\tRic$ is $\h$-equivariant, $\P_1(\h)$ is either trivial or isomorphic to $\Real^n$. The spaces $\P(\h)$ for irreducible Riemannian holonomy algebra $\h\subset\so(n)$ are computed recently in \cite{GalP}. In particular,
$\P_0(\h)\neq 0$ and $\P_1(\h)=0$ exactly for the holonomy algebras 
$\su(\frac{n}{2})$, $\sp(\frac{n}{4})$, $\spin(7)$ and $G_2$. Next, $\P_1(\h)\simeq\Real^n$ and $\P_0(\h)\neq 0$ exactly for the holonomy algebras $\so(n)$, $\u(\frac{n}{2})$ and $\sp(\frac{n}{4})\oplus\sp(1)$. For the rest of the Riemannian holonomy algebras (i.e. for the symmetric Berger algebras) it holds $\P_1(\h)\simeq\Real^n$ and $\P_0(\h)=0$.

To make the exposition complete, we give a description of the space $\R(\so(1,n+1))$ that follows immediately from
\cite{Al}. The space $\R(\so(1,n+1))$ admits the decomposition \eqref{razlR}. The complexification
$\R_0(\so(1,n+1))\otimes\Co$ of $\R_0(\so(1,n+1))$ is isomorphic to the $\so(n+2,\Co)$-module $V_{2\pi_2}$ (if $n\geq 3$).
Similarly, $\R_0(\so(1,3))\otimes\Co\simeq V_{4\pi_1}\oplus V_{4\pi'_1}$, $\R_0(\so(1,2))=0$, and $\R_0(\so(1,1))=0$.  It
holds $\R'(\so(1,n+1))\simeq(\odot^2\Real^{1,n+1})_0=\odot^2\Real^{1,n+1}/\Real\eta$, any $R\in \R'(\so(1,n+1))$ is of the
form $R=R_S$, where $S:\Real^{1,n+1}\to\Real^{1,n+1}$ is a symmetric linear map with zero trace and $$R_S(u,v)=Su\wedge v+
u\wedge Sv.$$ Similarly, $\R_1(\so(1,n+1))$ is spanned by the element $R=R_{\frac{\id}{2}}$, i.e. $R(u,v)=u\wedge v$. This
shows that an $(n+2)$-dimensional Lorentzian manifold ($n\geq 2$) with the holonomy algebra $\so(1,n+1)$ may be either
vacuum Einstein, or Einstein and not vacuum Einstein, or not Einstein. Finally, there are no vacuum Einstein indecomposable Lorentzian
manifolds of dimensions 2 or 3.

\section{Applications to  Einstein  and vacuum Einstein Lorentzian manifolds}\label{Einstein}

Now we are able to find all holonomy algebras of Einstein  and vacuum Einstein Lorentzian manifolds.

Let $R\in\R(\g^{1,\h})$ be as in Theorem \ref{thR}, then its Ricci  tensor $\Ric=\Ric(R)$ satisfies
\begin{align}\label{Ric1} \Ric(p,q)=&-\lambda,\quad \Ric(x,y)=\Ric(R_0)(x,y),\\
\label{Ric2} \Ric(x,q)=&\eta(x,\tRic(P)-v),\quad \Ric(q,q)=\tr T, \end{align} where $x,y\in\Real^n$ (recall that
$\Ric(u,v)=\tr(z\mapsto R(u,z)v)$).

Obviously, these equations imply that there is no three-dimensional indecomposable Einstein Lorentzian manifolds with
holonomy algebras contained in $\simil(1)=\Real\zr\Real$. Thus we may assume that  $n\geq 2$.

Here is a result from \cite{ESI}.

\begin{theorem}\cite{ESI}\label{thESI}
Let $(M,g)$ be a locally indecomposable Lorentzian Einstein manifold admitting a parallel distribution of isotropic lines.
Then the holonomy of $(M,g)$ is either of type 1 or 2. If the cosmological constant of $(M,g)$ is non-zero, then the
holonomy algebra of $(M,g)$ is of type 1. If $(M,g)$ admits locally a parallel isotropic vector field, then $(M,g)$ is
vacuum Einstein.
\end{theorem}

The  classification completes the following two theorems. 

\begin{theorem}\label{thRic0}
Let $(M,g)$ be a locally indecomposable $n+2$-dimensional Lorentzian manifold admitting a parallel distribution of
isotropic lines. If $(M,g)$ is vacuum Einstein, then one of the following holds:
\begin{description} \item[(1)] The holonomy algebra $\g$ of $(M,g)$ is of type 1, and in the decomposition \eqref{LM0B}
of the orthogonal part $\h\subset\so(n)$ at least one subalgebra $\h_i\subset\so(n_i)$ coincides with one of the Lie
algebras $\so(n_i)$, $\u(\frac{n_i}{2})$, $\sp(\frac{n_i}{4})\oplus\sp(1)$ or with a  symmetric
 Berger algebra.
\item[(2)] The holonomy algebra $\g$ of $(M,g)$ is of type 2, and in the decomposition \eqref{LM0B}
of the orthogonal part $\h\subset\so(n)$ each subalgebra $\h_i\subset\so(n_i)$ coincides with one of the Lie algebras
$\so(n_i)$, $\su(\frac{n_i}{2})$, $\sp(\frac{n_i}{4})$, $G_2\subset\so(7)$, $\spin(7)\subset\so(8)$. \end{description}
\end{theorem}

\begin{theorem}\label{thEinstein}
Let $(M,g)$ be a locally indecomposable $n+2$-dimensional Lorentzian manifold admitting a parallel distribution of
isotropic lines. If $(M,g)$ is Einstein and not vacuum Einstein, then the holonomy algebra $\g$ of $(M,g)$ is of type 1, and in
the decomposition \eqref{LM0B} of the orthogonal part $\h\subset\so(n)$ each subalgebras $\h_i\subset\so(n_i)$ coincides
with one of the Lie algebras $\so(n_i)$, $\u(\frac{n_i}{2})$, $\sp(\frac{n_i}{4})\oplus\sp(1)$ or with a symmetric
 Berger algebra. Moreover, $\h\subset\so(n)$ does not  annihilate any proper subspace of
$\Real^n$, i.e. in the decomposition \eqref{LM0A} it holds $n_{s+1}=0$.
\end{theorem}

Recall that the list of irreducible symmetric Berger algebras 
$\h\subset\so(n)$ can be obtained from the list of the holonomy algebras of irreducible Riemannian symmetric spaces (this list is given e.g. in \cite{Besse}) omitting $\so(n)$, $\u(\frac{n}{2})$ and $\sp(\frac{n}{4})\oplus\sp(1)$.

{\bf Proof of Theorems \ref{thRic0} and \ref{thEinstein}.} Fix a point $x\in M$ and let $\g$ be the holonomy algebra of
$(M,g)$ at this point. We identify $(T_xM,g_x)$ with $(\Real^{1,n+1},\eta)$. By the Ambrose-Singer Theorem \cite{Besse},
$\g$ is spanned by the images of the elements $R_\gamma=\tau_\gamma^{-1}\circ
R_y(\tau_\gamma\cdot,\tau_\gamma\cdot)\circ\tau_\gamma,$ where $\gamma$ is a piecewise smooth curve in $M$ starting at the
point $x$ and with an end-point $y\in M$, and $\tau_\gamma:T_xM\to T_yM$ is the parallel transport along $\gamma$. All
these elements belong to the space $\R(\g)$ and they can be given as in Theorem \ref{thR}. Suppose that $(M,g)$ is vacuum Einstein.
By Theorem \ref{thESI}, $\g$ is either of type 1 or 2. If $\g$ is  of type 2, then from \eqref{Ric1} and \eqref{Ric2} it
follows that each $R_\gamma$ satisfies $\lambda=0$, $v=0$, $\Ric(R_0)=0$, $\tRic(P)=0$, and $\tr T=0$. Hence the
orthogonal part $\h\subset\so(n)$ of $\g$ is spanned by the images of elements of $\R_0(\h)$ and $\P_0(\h)$. 
From \cite{GalP} it follows that $\h\subset\so(n)$ is spanned by the images of elements of $\R_0(\h)$.
Thus $\h$ is
the holonomy algebra of a vacuum Einstein Riemannian manifold. If  $\g$ is of type 1, then each $R_\gamma$ satisfies
$\lambda=0$, $v=\tRic(P)$, $\Ric(R_0)=0$, and $\tr T=0$. Hence for some element $R_\gamma$ it holds $\tRic(P)\neq 0$, i.e.
at list for one $\h_i\subset\so(n_i)$ in the decomposition \eqref{LM0B} it holds $\P_1(\h_i)\neq 0$.
 If $(M,g)$ is Einstein with the cosmological constant $\Lambda\neq 0$, then by Theorem \ref{thESI}, $\g$ is of type 1.
We get that  the curvature tensor $R_x$ at the point $x$ given by \eqref{thR} satisfies $\lambda=-\Lambda$, $v=\tRic(P)$,
and $\Ric(R_0)=\Lambda \eta|_{\Real^n\otimes\Real^n}$. Hence for each  $\h_i\subset\so(n_i)$ in the decomposition
\eqref{LM0B} it holds $\R_1(\h_i)\neq 0$, and $n_{s+1}=0$. $\Box$

{\bf Remark.} A simple version of Theorem \ref{thRic0} for Lorentzian manifolds with holonomy algebras of type 2 is proved
in \cite{ESI}, where the possibility for $\h_i\subset\so(n_i)$ to coincide with the holonomy algebra of a symmetric
Riemannian non-K\"ahlerian  space was not excluded.

\section{Examples of Einstein and vacuum Einstein Lorentzian metrics}\label{secExamples}

On  an $n+2$-dimensional Lorentzian manifold $(M,g)$ admitting a parallel distribution of isotropic lines there exist
local coordinates (the Walker coordinates) $x^0,...,x^{n+1}$ such that the metric $g$ has the form
\begin{multline}\label{Walker}
g=2dx^0dx^{n+1}+\sum^{n}_{i,j=1}h_{ij}(x^1,...,x^{n+1})dx^idx^j\\+2\sum^{n}_{i=1}u^i(x^1,...,x^{n+1}) dx^i
dx^{n+1}+f(x^0,...,x^{n+1})(dx^{n+1})^2,\end{multline} where $h(x^{n+1})=\sum^{n}_{i,j=1}h_{ij}(x^1,...,x^{n+1})dx^idx^j$
is a family of Riemannian metrics depending on the coordinate $x^{n+1}$ \cite{Walker}. The parallel distribution of
isotropic lines is defined by the vector field $\p_0$ (we denote $\frac{\p}{\p x^a}$ by $\p_a$).

In \cite{G-P} the Einstein equations for the general metric \eqref{Walker} are written down and some solutions of this
system under some assumptions on the coefficients are found. Of course, it is not possible to solve such system in
general.

In this section we show the existence of a local metric for each holonomy algebra obtained in the previous section. It is
easy to see that if the metric \eqref{Walker} is Einstein, then each Riemannian metric in the family $h(x^{n+1})$ is
Einstein with the same cosmological constant.

We will take some special $h_{ij}$ and $u^i$ in \eqref{Walker}, then the condition on the metric $g$ to be Einstein  will be equivalent to a system of equations on the function $f$. In each case we will show the existence of a proper function $f$ satisfying these equations. 
This will imply the existence of an Einstein metric with each holonomy algebra obtained above.

First consider the metric \eqref{Walker} with $h_{ij}$ independent of $x^{n+1}$ and $u^i(x^1,...,x^{n+1})=0$ for all
$i=1,...,n$. The holonomy algebras of such metrics are found in \cite{BB-I}. If  $f$ is sufficiently general, e.g. its Hessian is non-zero, then the holonomy algebra $\g$ of this metric is weakly-irreducible. If $\p_0 f=0$, then $\g=\h\zr\Real^n$, where $\h$ is
the holonomy algebra of the Riemannian metric $h$;  if $\p_0 f\neq 0$, then $\g=(\Real\oplus\h)\zr\Real^n$. In addition we
need to choose $h$ and $f$ to make $g$ Einstein or vacuum Einstein. The Ricci tensor $\Ric(g)$ for such metric has  the
following non-zero components:
\begin{align}
\Ric_{0\,n+1}&=\frac{1}{2}(\p_0)^2f,\\ \Ric_{i\,j}&=\Ric_{i\,j}(h),\quad i,j=1,...,n,\\
\Ric_{i\,n+1}&=\frac{1}{2}\p_0\p_if,\quad i=1,...,n,\\
\Ric_{n+1\,n+1}&=\frac{1}{2}\Big(f(\p_0)^2f-\Delta f\Big),\end{align} where $\Delta f= \sum_{i,j=1}^n h^{ij}\left(\p_{i}\p_{j}f-\sum_{k=1}^n \Gamma^k_{ij}\p_k f\right)$ is the Laplace-Beltrami operator of the metric $h$ applied to $f$. Suppose that $g$ is vacuum Einstein, then
the metric $h$ should be vacuum Einstein as well. Next, $\p_0f=0$ and $\Delta f=0$. 
Let $f$ be any harmonic function with non-zero Hessian, locally such functions always exist. E.g. if $h$ is flat, then we may take $f=(x^1)^2+\cdots+(x^{n-1})^2-(n-1)(x^n)^2$. Thus, the metric $g$ is vacuum Einstein and $\g=\h\zr\Real^n$.

Suppose that $g$ is Einstein with the cosmological constant $\Lambda\neq 0$, then  $h$ is Einstein with the same
cosmological constant $\Lambda$. Next, $\frac{1}{2}(\p_0)^2f=\Lambda$, $\p_0\p_if=0$ and
$\frac{1}{2}\big(f(\p_0)^2f-\Delta f)=\Lambda f$. We get that $f=\Lambda\, (x^0)^2+x^0 f_1(x^{n+1}) +f_0(x^1,...,x^{n+1})$ such
that $\Delta f_0=0$. Taking $f_0$ to be harmonic with non-zero Hessian, and $f_1=0$, 
we get that the metric $g$ is
Einstein with the cosmological constant $\Lambda$ and $\g=(\Real\oplus\h)\zr\Real^n$.

It is impossible to construct in this way vacuum Einstein metrics with the holonomy algebras of type~1.

 In \cite{Gal5} for each weakly-irreducible not irreducible holonomy algebra is constructed
a metric of the form \eqref{Walker} with $h_{ij}(x^1,...,x^{n+1})=\delta_{ij}$, i.e. each  Riemannian metric in the family
$h(x^{n+1})$ is flat. The Ricci tensor $\Ric(g)$ for such metric has the following components:
\begin{align}\label{Ric9}
\Ric_{0\,n+1}=&\frac{1}{2}(\p_0)^2f,\\ \Ric_{i\,j}=&0,\quad i,j=1,...,n,\\\label{Ric11}
\Ric_{i\,n+1}=&\frac{1}{2}\Big(\p_0\p_if-\sum_{j=1}^n\p_j(\p_ju^i-\p_iu^j)\Big),\quad i=1,...,n,\\
\Ric_{n+1\,n+1}=&\frac{1}{2}\Big(\big(f-\sum_{i=1}^n
(u^i)^2\big)(\p_0)^2f-\sum_{i=1}^n(\p_i)^2f+2\sum_{i=1}^n\p_i\p_{n+1}u^i\\\nonumber&+\sum_{i,j=1}^n(\p_ju^i-\p_iu^j)^2+(\p_0f)\sum_{i=1}^n\p_iu^i+2\sum_{i=1}^nu^i\p_0\p_if\Big).\end{align}
Now we recall the algorithm of the construction from \cite{Gal5}.

Let $\h\subset\so(n)$ be the holonomy algebra of a Riemannian manifold. We get the decompositions \eqref{LM0A} and
\eqref{LM0B}. We will assume that the basis $e_1,...,e_n$ of $\Real^n$ is compatible  with the decomposition of $\Real^n$.
Let $n_0=n_1+\cdots+n_s=n-n_{s+1}$. We see that $\h\subset\so(n_0)$ and $\h$ does not annihilate any proper subspace of
$\Real^{n_0}$. We will always assume that  the indices $i,j,k$ run from $1$ to $n$, the indices $\i,\j,\k$ run from $1$ to
$n_0$, and the indices $\ii,\jj,\kk$ run from $n_0+1$ to $n$. We will use the
Einstein rule for sums.

Let $(P_\alpha)_{\alpha=1}^N$ be linearly independent elements of $\P(\h)$ such that the subset
$\{P_\alpha(u)|1\leq\alpha\leq N,\,u\in \Real^{n}\}\subset\h$ generates the Lie algebra $\h$. For example, it can be any
basis of the vector space $\P(\h)$.  For each $P_\alpha$ define the numbers $P_{\alpha \j \i}^\k$ such that
$P_\alpha(e_\i)e_\j=P_{\alpha \j \i}^\k e_\k.$ Since $P_\alpha\in\P(\h)$, we have
\begin{equation}\label{LM1} P_{\alpha\k\i}^\j=-P_{\alpha\j\i}^\k \text{ and }
 P_{\alpha\j\i}^\k+P_{\alpha\k\j}^\i+P_{\alpha\i\k}^\j=0.\end{equation}
 It holds $\tRic(P_\alpha)=\tRic(P_\alpha)^\k e_\k$, where $\tRic(P_\alpha)^\k=\sum_\i P_{\alpha\i\i}^\k$.
Define the following numbers
\begin{equation}\label{LM2A} a_{\alpha \j \i}^\k=\frac{1}{3\cdot(\alpha-1)!}\left(P_{\alpha \j \i}^\k+P_{\alpha \i \j}^\k\right).
\end{equation}
We have \begin{equation}\label{LM3} a_{\alpha \j \i}^\k=a_{\alpha \i \j}^\k.\end{equation} From \eqref{LM1} it follows
that
\begin{equation}\label{LM4} P_{\alpha \j \i}^\k=(\alpha-1)!\left(a_{\alpha \j \i}^\k-a_{\alpha \k \i}^\j\right) \text{ and }
 a_{\alpha \j \i}^\k+a_{\alpha \k \j}^\i+a_{\alpha \i \k}^\j=0.\end{equation}
Define the functions  \begin{equation}\label{LM6A}u^\i=a_{\alpha \j \k}^\i x^\j x^\k(x^{n+1})^{\alpha-1}\end{equation} and
set $u^\ii=0$. We choose the function $f$ to make the holonomy algebra $\g$ of the metric $g$ to be weakly-irreducible.
 If $\p_0 f=0$, then $\g$ is either of type 2 or type 4; if
$\p_0 f\neq 0$, then $\g$ is either of type 1 or type 3. We will make $g$ to be vacuum Einstein, then $\g$ will be either of
type 2 or type 1, i.e. it will equal either to $\h\zr\Real^n$ or to $(\Real\oplus\h)\zr\Real^n$.

Note that
\begin{equation}\label{svojstvou} \p_\j u^\i-\p_\i u^\j=\frac{2}{(\alpha-1)!}P^\i_{\alpha\j\k}x^\k(x^{n+1})^{\alpha-1},
\qquad\p_\i u^\i=-\frac{2}{3((\alpha-1)!)}\sum_{\k}\tRic(P_\alpha)^\k x^\k(x^{n+1})^{\alpha-1}.\end{equation}

Suppose that $g$ is vacuum Einstein. Then, first of all, the equality \eqref{Ric9}  implies that
$f=x^{0}f_1+f_0$, where $f_0$ and $f_1$ are functions of $x^1,...,x^{n+1}$.

First suppose that $\g$ is of type 2, i.e. $\p_0f=0$ and $f_1=0$. Substituting this and \eqref{svojstvou} into the
equation $\Ric=0$, we get the following equations $$\sum_{\j} P^\i_{\alpha \j\j}=0,\qquad \sum_{i}(\p_i)^2
f_0=\sum_{\i,\j}\Big(\frac{2}{(\alpha-1)!}P^\i_{\alpha \j\k} x^\k (x^{n+1})^{\alpha-1}\Big)^2.$$ The first equation is
equivalent to the condition $\tRic(P_\alpha)=0$. Clearly, the function \begin{equation}\label{f}
f_0=\frac{1}{3}\sum_{\i,\j}\Big(\frac{1}{(\alpha-1)!}P^\i_{\alpha \j\k} (x^\k)^2 (x^{n+1})^{\alpha-1}\Big)^2\end{equation}
satisfies the second equation. In order to make $\g$ weakly-irreducible we add to the obtained $f_0$ the harmonic function
$(x^1)^2+\cdots+(x^{n-1})^2-(n-1)(x^{n})^2$ (it is not necessary to do this if $n_0=n$).

 Thus we get a new example of the vacuum Einstein metric with the holonomy algebra
$\h\zr\Real^n$, where $\h$ is the (not necessary irreducible) holonomy algebra of a vacuum Einstein Riemannian manifold.

Suppose now that $\g$ is of type 1. The equation $\Ric=0$ is equivalent to the following system of equations:
\begin{align*} \p_\i f_1=&\frac{2}{(\alpha-1)!}\tRic(P_\alpha)^\i(x^{n+1})^{\alpha-1},\qquad \p_\ii f_1=0,\\
\sum_i(\p_i)^2 f=&\sum_{\i,\j}\Big(\frac{2}{(\alpha-1)!}P^\i_{\alpha \j\k} x^\k (x^{n+1})^{\alpha-1}\Big)^2-
\sum_\i\frac{4}{3((\alpha-2)!)}\tRic(P_\alpha)^\i x^\i(x^{n+1})^{\alpha-2}\\&-
f_1\sum_\i\frac{2}{3((\alpha-1)!)}\tRic(P_\alpha)^\i x^\i(x^{n+1})^{\alpha-1}+2u^\i\p_\i f_1 .\end{align*} We may take
$f_1=\sum_\i\frac{2}{(\alpha-1)!}\tRic(P_\alpha)^\i x^\i(x^{n+1})^{\alpha-1}$. Substituting it into the last equation, we
obtain the  equation of the form $$\sum_i(\p_i)^2 f_0=G,$$ where $G$ is a polynomial of the variables $x^\i$ and
$x^{n+1}$, and it is of degree at most 2 in the variables $x^\i$ and of degree at most  $2(N-1)$ in $x^{n+1}$. Let
$f_0=\sum_{\beta=0}^{2(N-1)}f_{0\beta}\cdot(x^{n+1})^\beta$, then each $f_{0\beta}$ satisfies the equation $$\sum_\i
(\p_\i)^2 f_{0\beta} =G_\beta,$$ where $G_\beta$ is a polynomial of $x^\i$ of degree at most 2. Thus we need to find
solutions of a number of the Poisson equations \begin{equation}\label{FH}\sum_\i (\p_\i)^2 F=H,\end{equation} where $H$ is
a polynomial of $x^\i$ of degree at most 2. Let us show a simple way to find a polynomial solution of such equation. Let
$H_1=H-\frac{1}{2}(x^\i)^2(\p_\i)^2H,$ then $H=H_1+\frac{1}{2}(x^\i)^2(\p_\i)^2H$ and $(\p_\i)^2H_1=0$ for each $\i$.
Next, let $H_2=H_1-x^1\p_1H_1$, then $H_1=H_2+x^1\p_1H_1$, $\p_1 H_2=0$, and $(\p_\i)^2H_1=0$ for each $\i$. Now it is
obvious that the function $$F=\frac{1}{2}(x^1)^2 H_2+\frac{1}{6}(x^1)^3 \p_1 H_1+\frac{1}{24}(x^\i)^4(\p_\i)^2H$$ is a
solution of the equation \eqref{FH}. In order to make $\g$ weakly-irreducible we add to the obtained $f_0$ the harmonic
function $(x^1)^2+\cdots+(x^{n-1})^2-(n-1)(x^{n})^2$.

Thus we get an example of the vacuum Einstein metric with the holonomy algebra $(\Real\oplus\h)\zr\Real^n$, where $\h$ is the
(not necessary irreducible) holonomy algebra of a Riemannian manifold such that $\P_1(\h)\neq 0$, in other words, in the
decomposition \eqref{LM0B} at least one $\h_i$ is the holonomy algebra of a not vacuum Einstein Riemannian manifold.

We have proved the following theorem.

\begin{theorem} Let $\g$ be any algebra as in Theorem \ref{thRic0} or \ref{thEinstein}, then there exists an $n+2$-dimensional
Einstein (resp., vacuum Einstein) Lorentzian manifold with the holonomy algebra $\g$. \end{theorem}

\begin{ex} In \cite{Gal5,ESI} we constructed metrics with the holonomy algebras $\g^{2,G_2}\subset\so(1,8)$ and
$\g^{2,\spin(7)}\subset\so(1,9)$. In these constructions $N=1$ and $f=0$. Choosing in these constructions
$f=\frac{1}{3}\sum_{i,j=1}^{n}(P^i_{j k} (x^k)^2)^2$,  we obtain vacuum Einstein metrics with the holonomy algebras
$\g^{2,G_2}\subset\so(1,8)$ and $\g^{2,\spin(7)}\subset\so(1,9)$. Similarly, using the metric constructed in \cite{Gal5} on page 1033,
it is easy to construct vacuum Einstein metric with the holonomy algebra $\g^{1,\rho(\so(3))}\subset\so(1,6)$, where
$\rho:\so(3)\to\so(5)$ is the irreducible representation of $\so(3)$ on $\Real^5$. For this it is enough to choose $f=\left(-\frac{20}{3}x^1-2x^3\right)x^0+f_0$, where $f_0$ is a polynomial of degree 4 that can be easily found using the above algorithm.
Note that $\rho(\so(3))\subset\so(5)$ is a symmetric Berger algebra and it is the holonomy algebra of the symmetric space ${\rm SL}(3,\Real)/{\rm SO}(3,\Real)$. \end{ex}

\section{Lorentzian manifolds with totally isotropic Ricci operator} 

Let $R\in\R(\g^{1,\h})$ be as in Theorem \ref{thR}. Consider its Ricci operator $\Ric=\Ric(R):\Real^{1,n+1}\to\Real^{1,n+1}$ defined by
$$\eta(\Ric (x),y)=\Ric(x,y),$$ where $x,y\in\Real^{1,n+1}$ and on the right hand side $\Ric$ denotes the Ricci tensor of $R$.
It is easy to check that
\begin{align}\label{Ric3} \Ric(p)=&-\lambda p,\quad \Ric(x)=\eta(x,\tRic(P)-v)p+\Ric(R_0)(x),\\
\label{Ric4} \Ric(q)=&(\tr T)p+\tRic(P)-v-\lambda q,\end{align} where $x\in\Real^n$. 

A Lorentzian manifold $(M,g)$ is called {\it totally Ricci-isotropic} if the image of its Ricci operator is isotropic, equivalently,
$\eta(\Ric(x),\Ric(y))=0$ for all $x,y\in \Real^{1,n+1}$. Obviously, any vacuum Einstein Lorentzian manifold is totally Ricci-isotropic.
If $(M,g)$ is a spin manifold and it admits a parallel spinor, then it is totally Ricci-isotropic (but not necessary vacuum Einstein, unlike in the Riemannian case) \cite{Bryant,FF00}.

\begin{theorem}\label{thRicisotr}
Let $(M,g)$ be a locally indecomposable $n+2$-dimensional Lorentzian manifold admitting a parallel distribution of
isotropic lines. If $(M,g)$ is totally Ricci-isotropic, then its holonomy algebra is the same  as in Theorem \ref{thRic0}. 
\end{theorem}

{\it The proof} of this theorem is similar to the proofs of the above Theorem \ref{thESI} from \cite{ESI} and Theorems \ref{thRic0} and \ref{thEinstein}. $\Box$

Using results of Section \ref{curv} it is easy to prove the following theorem.

\begin{theorem}\label{thRicisotr1} Let $(M,g)$ be a locally indecomposable $n+2$-dimensional Lorentzian manifold admitting a parallel distribution of
isotropic lines. If the holonomy algebra of $(M,g)$ is of type 2 and in the decomposition \eqref{LM0B}
of the orthogonal part $\h\subset\so(n)$ each subalgebra $\h_i\subset\so(n_i)$ coincides with one of the Lie algebras
 $\su(\frac{n_i}{2})$, $\sp(\frac{n_i}{4})$, $G_2\subset\so(7)$, $\spin(7)\subset\so(8)$, then $(M,g)$ is totally Ricci-isotropic.
\end{theorem}

Note that this theorem can  be also proved by the following argument. Locally $(M,g)$ admits a spin structure. From \cite{ESI,Leistner} it follows that any manifold $(M,g)$ with the holonomy algebra as in the theorem admits locally a parallel spinor, hence $(M,g)$ is totally Ricci-isotropic.


\end{document}